\documentclass[12pt,reqno]{amsart}
\usepackage{latexsym,amsmath,amsthm,amssymb,amsfonts}
\usepackage{setspace}
\usepackage{color}
\usepackage{tikz}
\numberwithin{equation}{section}
\newtheorem{theorem}{Theorem}[section]
\newtheorem{lem}[theorem]{Lemma}
\newtheorem{thm}[theorem]{Theorem}
\newtheorem{pro}[theorem]{Proposition}
\newtheorem{cor}[theorem]{Corollary}

\newtheorem{rem}[theorem]{Remark}

\def\S2{\mathbb{S}^2}
\def\s{\,\,\,\,}
\def\lan{\langle}
\def\ran{\rangle}

\def\mv{1.7ex}
\def\endproof{$\hfill\Box$\\}
\def\R{\mathbb{R}}
\def\C{\mathbb{C}}
\def\vol{\mu}

\title{\bf Metrics  On A Surface With  Bounded Total Curvature}
\author{Yuxiang Li, Jianxin Sun, Hongyan Tang}

\address{\newline
Yuxiang Li:
 Department of Mathematical Sciences, Tsinghua University, Beijing 100084, P.R. China.
{\tt Email:yxli@math.tsinghua.edu.cn}
\newline
\newline
Jianxin Sun:
 Academy of Mathematics and Systems Science, Chinese Academy of Sciences, Beijing 100190, P.R. China.
{\tt Email:sunjianxin201008@126.com}
\newline
\newline
Hongyan Tang:
Department of Mathematical Sciences, Tsinghua University, Beijing 100084, P.R. China.
{\tt Email:hytang@math.tsinghua.edu.cn}}

\date{}
\begin{document}
\maketitle

\begin{abstract}
Let $g=e^{2u}(dx^2+dy^2)$ be a conformal metric defined on the unit disk of $\C$. We   give an estimate of $\|\nabla u\|_{L^{2,\infty}(D_\frac{1}{2})}$ when $\|K(g)\|_{L^1}$
is small and $\frac{\vol(B_r^g(z),g)}{\pi r^2}<\Lambda$
  for any $r$ and $z\in D_\frac{3}{4}$. Then we will use this estimate
to study the Gromov-Hausdorff convergence of a conformal metric sequence with bounded $\|K\|_{L^1}$ and
give some applications.
\end{abstract}

\section{Introduction}
In \cite{Shioya}, Takashi Shioya studied the convergence of surfaces  and proved the following.
\begin{thm}[Shioya]
Let $\mathcal{M}(C,d)$ denote the set of isometry classes of closed 2-dimensional Riemannian manifolds  with diameter $\leq d$ and total absolute curvature $\|K\|_{L^1}\leq C$. For any positive $C$ and $d$,
$\mathcal{M}(C,d)$ is precompact with respect to the Gromov-Hausdorff distance. 
\end{thm}
In this paper, 
we use the theory of moduli space and Brezis-Merle 
Theorem to give more information of the convergence in case of Riemann 
surfaces. The idea is the following. By the theory of moduli space, we can divide a Riemann surface into finitely many
parts, and each part is conformal to a disk or $D\setminus D_{r_k}$ with $r_k\rightarrow 0$, where $D$ is the unit disk, and $D_{r_k}=D_{r_k}(0)$. So if we get the convergence
of $D$ and $D\setminus D_{r_k}$, we can study the convergence of
a surface sequence by a standard blowup analysis.  In addition, the converging behavior
of a metric sequence on $D\setminus D_{r_k}$ can be also obtained from
a precise bubble tree argument of a metric sequence defined on the
unit disk. So the key point to such a problem is to understand the
detail of a conformal metric sequence defined on a disk. 
The main topic of this paper   is the  
first step of this program.
We need to point out that similar problems under different assumptions have already been
studied by many authors \cite{Cheng-Lin1,Cheng-Lin2, Chen, Li-Shafrir, Zhang}.

Set $g_k=e^{2u_k}g_{euc}$, where $g_{euc}=dx^2+dy^2$.
Let $K(g_k)$ be the Gauss curvature of $g_k$. We have the following Gauss curvature equation:
$$
-\Delta u_k=K(g_k)e^{2u_k}.
$$
By Brezis-Merle's result (see section 2), we can decompose $u_k$ into a sum of a harmonic function $w_k$ and a function $v_k$ which is bounded in
$W^{1,p}$.  Thus, to understand the convergence of $u_k$, we
only need to study the behavior of the mean value of $w_k$. 

The first result in this direction was obtained by H\'elein. In \cite{Helein}, H\'elein proved that
 a conformal immersion sequence from $D$ into $\R^n$ with induced metrics
$g_k = e^{2u_k}g_{euc}$, with  small $L^2$-norm of second fundamental and area bounds, is bounded in $W^{2,2}_{loc}(D,\R^n)$ 
and converges to a conformal immersion or a point. Moreover, if
the limit is not a point, then $u_k$ is bounded in $L^\infty_{loc}(D)$. 
In H\'elein's case, $u_k$ can be decomposed into the sum of
a $L^\infty$  function and a harmonic function, so he was able to
discuss the $W^{1,2}$ convergence and $L^\infty$ bound of $u_k$.

We call the sequence collapse if the limit is a point. Usually, collapse is 
the most difficult part in the analysis.
When it happens, one will be lucky to find a constant $c_k$, such that $u_k-c_k$ converges.
Unfortunately, we have the following counterexample. Let $f_k=\frac{e^{kz}}{a_k}:D\rightarrow\C$, and $g_k=e^{2u_k}g_{euc}$ is the induced metric.
We choose $a_k$  such that 
$
\vol(D,g_k)=\int_D\frac{e^{2kx}}{a_k^2}=1.
$
Then for any $c_k$, $u_k-c_k$ can not converge. However, in a recent result \cite{Li-Wei-Zhou}, we showed
that this is almost the only counterexample.  Specifically, we proved that if $\|A_k\|_{L^2}\rightarrow 0$ and $\|\nabla u_k\|$ is not bounded, then we can find $x_k$ and $r_k$, such that $u_k(x_k+r_kx)$ converges to a linear function. As an application, we proved that when $\frac{\vol(f_k(D)\cap B_r^n(p))}{\pi r^2}<C$
for any $B_r^n(p)$, we can find $c_k$ such that $u_k-c_k$
converges weakly in $W^{1,2}$.

We will show that a similar argument works in the intrinsic case.
Precisely, we will use  the John-Nirenberg radius defined in
\cite{Li-Wei-Zhou} to measure the Harnack properties of a function satisfying
an elliptic equation. The main observation is, when $\Delta u_k$ is very
closed to 0 in $L^1$-norm, then we can find sequence $u_k(x_k+r_k x)$
which converges to a linear function. 

The first result of this paper is   the following.

\begin{thm}\label{main1}
Let $g=e^{2u}g_{euc}$ be a smooth metric defined on $D$.
We assume
\begin{equation}\label{volume.assumption}
\frac{\vol(B_r^{g}(z),g)}{\pi r^2}<\Lambda_1,\s \forall r>0\s and \s z\in D_\frac{3}{4}.
\end{equation}
Then there exists an $\epsilon_0$, such that if 
\begin{equation}\label{epsilon}
\int_D|K(g)|d\mu_{g}\leq\epsilon_0<\frac{4}{3}\pi,
\end{equation}
then
$$
\|\nabla u\|_{L^{2,\infty}(D_\frac{1}{2})}<C.
$$
Moreover, we have
$$
\|u\|_{W^{1,p}(D_\frac{1}{2})}\leq C
(1+|\log \vol(D_\frac{1}{2},g)|),
$$
where $p\in(1,2)$.
\end{thm}

\begin{rem}
By Corollary \ref{BH2} and Poincar\'e inequality, \eqref{volume.assumption}
and \eqref{epsilon} imply that 
$$
\int_{D_\frac{1}{2}}e^{q(u-c)}<C(q,\epsilon_0,\Lambda_1),
$$
where $c$ is the mean value of $u$ along $\partial D_\frac{1}{4}$,
and $q<\frac{4\pi}{\epsilon_0}$.
\end{rem}

As a consequence, we have the following corollary.
\begin{cor}\label{main2}
Let $g_k=e^{2u_k}g_{euc}$ be a smooth metric defined on $D$.
We assume $(D, g_k)$ can be extended to a complete surface with 
total Gauss curvature uniformly bounded. If 

$$
\int_D|K(g_k)|d\mu_{g_k}<\tau\leq\epsilon_0,
$$
then all  results in Theorem \ref{main1} hold.
\end{cor}

Together with the theory of moduli space of Riemann surface, this corollary may help us to study the convergence of 
a Riemann surface sequence with bounded total Gauss curvature.

It is easy to check that if $\vol(D, g_k)$ is uniformly bounded, then $c_k$ must be uniformly bounded above. If $\vol(D, g_k)\rightarrow 0$, then $c_k\rightarrow-\infty$, hence $\|e^{u_k}\|_{W^{1,p}(D_r)}\rightarrow 0$. Then by the trace embedding inequality, we get the following corollary.
\begin{cor}\label{collapse}
 Let $g_k$ be as in Corollary \ref{main2} and assume $\vol(D, g_k)<\Lambda_2$ in addition. Then, after passing to a subsequence, $e^{u_k}$ converges in $L^q(D_r)$ for any $q<\frac{4\pi}{\tau}$, and for any $r\in(0,1)$ and $p\in(1,2)$, one of the following two alternatives
holds:
\begin{itemize}
\item[{\rm (a)}] $u_k$ converges weakly in $W^{1,p}(D_r)$;
\item[{\rm (b)}]  $e^{u_k}$ converges to 0 in $L^q(D_r)$ and both  
$\vol(g_k,D_r)$   and
 $diam(D_r, g_k)$ converge to 0.
\end{itemize}
  Moreover, when $\rm(b)$ happens, we can find $c_k$, such that
 $u_k-c_k$ converges weakly in $W^{1,p}(D_r)$.

 \end{cor}

We need to study the limit distance space.
For this sake, we denote
$$
\mathcal{RF}^p(D)=\{u\in W^{1,p}(D):\Delta u \mbox{ is a finite Radon measure  }\}.
$$
Here $\Delta u$ is a Radon measure means that there exists 
a Radon measure $\vol$,
such that for any $\varphi\in\mathcal{D}$, 
$$
\int_D\nabla u\nabla\varphi=\int_D\varphi d\mu.
$$
Given $u\in\mathcal{RF}^p(D)$, we define
$$
d_u(x, y)=\inf_{\gamma\in P(x, y)}\int_\gamma e^{u},
$$
where $P(x, y)$ denotes the set of all piecewise smooth paths joining $x$ and $y$.
We will show that 
\begin{thm}\label{main3}
For any $u\in\mathcal{RF}^p(D)$ with $p\in(1,2)$,
$d_u$ is a metric. Moreover, if $g_k=e^{2u_k}g_{euc}$
is a smooth metric satisfying
$$
\int_D|K(g_k)|d\mu_{g_k}<\epsilon_0,
$$
and $u_k$ converges to $u$ weakly in $W^{1,p}$, then 
$(\overline{D_r},d_{g_k})$ converges to $(\overline{D_r},d_u)$ in the sense of Gromov-Hausdorff
distance.
\end{thm}

To make some preparations for
the blowup analysis for a surface sequence with 
finite total curvature in a forthcoming paper, we prove the following result.
 \begin{thm}\label{neck.estimate}
Let $g=e^{2u}g_{euc}$ on $D_4\setminus D_\frac{1}{4}$ with
$$
\|\nabla u\|_{L^{2,\infty}(D_4\setminus D_\frac{1}{4})}<\Lambda_3.
$$
Set $c$ to be the mean value
of $u$ on $\partial D_\frac{3}{2}$.
Then, there exists a constant $\epsilon_2
=\epsilon_2(\Lambda_3)$,
such that if
$$
\int_{D_4\setminus D_\frac{1}{2}}|K|e^{2u}<\epsilon_2,
$$
then
\begin{equation}\label{diam}
C_1\leq \frac{d_{g}(e^{i\theta},2e^{i\theta'})}{e^{c}}\leq C_2, \s\forall \theta,
\theta',
\end{equation}
and
$$
C_1'\leq \frac{\vol(D_2\setminus D, g)}{e^{2c}}\leq C_2', 
$$
where $C_1$, $C_1'$, $C_2$, $C_2'$ only depend on $\Lambda_3$
and $\epsilon_2$.
\end{thm}

In the last section, we will give some applications. First, we prove that under assumption of bounded total curvature, if a Riemann surface
sequence converges to a Riemann surface in the Gromov-Hausdorff distance, then the conformal  class and the volume converge.
Second, we will show that when $\|K(g_k)-1\|_{L^1}$ is small, 
$\{(S^2,g_k)\}$  is precompact  in the Gromov Hausdorff distance and in 
$W^{1,p}$ weakly. At the end, we prove that when $K(g_k)\geq 1$
and the $(S^2,g_k)$ does not collapse, then the Gromov-Hausdorff
limit is a topological sphere.

\section{Preliminary}

The arguments of this paper base on the following theorem.
\begin{thm}\label{BH}\cite{Brezis-Merle,Helein}\label{Briezis} Let $v$ be the solution of 
\begin{equation}\label{BH}
-\Delta v=f,\s v|_{\partial D}=0,
\end{equation}
where $f\in L^1$. Then for any $\epsilon>0$, we have
\begin{equation}\label{BH.estimate}
\|\nabla v\|_{L^{2,\infty}(D)}<C
\|f\|_{L^1(D)},\s and\s\int_\Omega e^{\frac{(4\pi-\epsilon)|v|}{\|f\|_{L^1(D)}}}<\frac{16\pi^2
}{\epsilon^2}.
\end{equation}
\end{thm}

\begin{rem}
Note that for any $p\in(1,2)$ and  $D_r(x)\subset D$,  
$$
\left(r^{2-p}\int_{D_r(x)}|\varphi|^p\right)^\frac{1}{p}<C\|\varphi\|_{L^{2,\infty}(D)}.
$$
Thus, if $v$ is the solution of \eqref{BH}, we have
$$
r^{2-p}\int_{D_r(x)}|\nabla v|^p dx\leq C\|f\|_{L^1},\s \forall D_r(x)\subset D.
$$
\end{rem}

Brezis-Merle's result can be considered as the $W^{2,1}$-version of
Trudinger inequality:
$$
\sup_{v\in W^{1,2}_0(\Omega),\|\Delta v\|_{L^1}=1}\int_\Omega
e^{(4\pi-\epsilon)|v|}<C(\epsilon).
$$
When $u\notin W^{1,2}_0$ with $|\Delta u|\in L^1$, 
we 
set $v\in W^{1,2}_0$ to be a solution of
$$
\Delta v=\Delta u.
$$
We have the decomposition:
$$
u=v+w,
$$
where  $w$ is harmonic and $v$ satisfies \eqref{BH.estimate}.
Then, by mean value theorem of harmonic function, we have the following:

\begin{cor}\label{BH2}
Let $u$ be a solution of 
\begin{equation*}
-\Delta u=f,\s x\in D.
\end{equation*}
If there exist $\epsilon_1$ and $\Lambda_3$ such that 
$$
\|f\|_{L^1(D)}<\epsilon_1,\s \|u\|_{L^1(D)}<\Lambda_3, 
$$
then
$$
\|e^{|u|}\|_{L^q(D_r)}<C(q,\Lambda_3,r).
$$
for any $q<\frac{4\pi}{\epsilon_1}$.
\end{cor}

For a function $\phi\in L^1(\Omega)$, we define $\rho$ as follows:
\begin{equation}\label{defi.JNR}
\rho(\phi, x_0,\Omega,\lambda)=\sup\{r:
D_t(x_0)\subset\Omega,\s \frac{1}{|D_t(x_0)|}\int_{D_t(x_0)}\left|\phi-
\frac{1}{|D_t(x)|}\int_{D_t(x_0)}\phi\right|dx\leq \lambda,\s \forall t<r\}.
\end{equation}
We define $\rho(\phi, x_0,\Omega,\lambda)=0$, if
$$
\varliminf_{r\rightarrow 0}\frac{1}{|D_r(x_0)|}\int_{D_r(x_0)}\left|\phi-
\frac{1}{|D_r(x_0)|}\phi\right|dx>\lambda.
$$

For a harmonic function, we have the following result.

\begin{lem}\label{gradient}
Let $u$ be harmonic on $D$ with $\rho(u, x,D,\lambda)>a>0$ for any
$x\in D_\frac{1}{2}$ and some positive number $a$. Then 
$$
\|\nabla u\|_{C^m(D_\frac{1}{2})}<C(a,\lambda, m).
$$
\end{lem}

\proof
Let $c=\frac{1}{|D_a(0)|}\int_{D_a(0)}udx$. It follows from 
$\rho(u,x,D,\lambda)>a$ that 
$$
\frac{1}{|D_a(0)|}\int_{D_a(0)}|u-c|dx\leq\lambda.
$$
By the mean value formula, we have
$$
\|u-c\|_{L^\infty(D_\frac{a}{2}(0))}<C.
$$
Since
$$
\frac{1}{|D_a(p)|}\int_{D_a(p)}|(u-c)-\frac{1}{|D_a(p)|}\int_{D_a(p)}(u-c)|dx\leq \lambda,\s\forall p\in D_\frac{a}{2}(0),
$$
by the mean value formula again,
we get
$$
\int_{D_a(p)}|u-c|dx\leq \int_{D_a(p)}|u-c-(u(p)-c)|dx+|u(p)-c||D_a(p)|.
$$
Hence $\|u-c\|_{L^\infty(D_{\frac{a}{2}2})}<C$. Then we are able to prove
$$
\|u-c\|_{L^\infty(D_{\frac{a}{2}k}(0))}<C
$$
whenever $\frac{a}{2}k<1$.
\endproof

Let $u=v+w$, where $v\in W^{1,2}_0(D)$ and $w$ is harmonic.
By Theorem \ref{BH} and Sobolev inequality, we have
\begin{equation}\label{BMO}
\frac{1}{|D_r(x)|}\int_{D_r(x)}\left|v-
\frac{1}{|D_r(x)|}\int_{D_r(x)}v\right|dx<Cr^{2-p}\int_{D_r(x)}|\nabla v|^pdx<C\|f\|_{L^1(D)}.
\end{equation}
Thus, if $\rho(u,x,D,\lambda)>a$, we get
$$
\rho(w,x,D,\lambda+C\|f\|_{L^1})>a.
$$
Then, we get the following corollary.

\begin{cor}
Suppose $u_k\in W^{1,p}(D)$ solve the equation
$$
-\Delta u_k=f_k
$$ 
in the sense of distribution. Assume 
$$\|f_k\|_{L^1(D)}<\epsilon_1,\s and\s
\inf_{x\in D_\frac{1}{2}}\rho(x,u_k,D,\lambda)>a>0,
$$
then  $\nabla u_k$ is bounded in $L^{2,\infty}(D_\frac{1}{2})$.
\end{cor}

\begin{lem}\label{linear.function}
Suppose $-\Delta u_k=f_k$, and assume
$$
\int_D|f_k|dx\rightarrow 0,\s
\inf_{x\in D_\frac{1}{2}}\rho(u_k,x,D,\lambda)\rightarrow 0.
$$
Then after passing to a subsequence, we can find $x_k\in D_\frac{1}{2}$, $r_k\rightarrow 0$,
and $c_k \in \mathbb{R}$,
such that $u_k(x_k+r_k x)-c_k$  converges to a nontrivial linear function weakly in 
$W^{1,p}(D_R)$ for any $R$ and $p\in(1,2)$.
\end{lem} 

\proof
Let $y_k\in D_\frac{1}{2}$, s.t.
$\rho(y_k,u_k,D,\lambda)\rightarrow 0$.
For simplicity, we denote $\rho(x,u_k,D,\lambda)$ by $\rho_k(x)$.

Put $x_k\in D_\frac{2}{3}$, such that
$$
\frac{\rho_k(x_k)}{2/3-|x_k|}=\inf_{x\in D_\frac{2}{3}}\frac{\rho_k(x)}{2/3-|x|}:=\tau_k.
$$
Noting that $$\tau_k\leq \frac{\rho_k(y_k)}{\frac{2}{3}-|y_k|}\rightarrow 0,$$
we have $\rho_k(x_k)=\tau_k(\frac{2}{3}-|x_k|)\rightarrow 0$, and hence for any fixed $R$
$$
D_{R\rho_k(x_k)}(x_k)\subset  D_{\frac{2}{3}-|x_{k}|}(x_{k})\subset D_\frac{2}{3},
$$
when $k$ is sufficiently large.
Then, for any $x\in D_{R\rho_k}(x_k)$ we have
\begin{eqnarray*}
\frac{\rho_k(x)}{\rho_k(x_k)}&\geq& \frac{\frac{2}{3}-|x|}{\frac{2}{3}-|x_k|}\geq \frac{\frac{2}{3}-|x_k|-|x-x_k|}{\frac{2}{3}-|x_k|}\\
&\geq&1-\frac{R\rho_k(x_k)}{\frac{2}{3}-|x_k|}\\
&=&1-R\tau_k.
\end{eqnarray*}
Hence, as $k$ is large enough, there holds
 $$\frac{\rho_k(x)}{\rho_k(x_k)}>\frac{1}{2}.$$

Set $r_k=\rho_k(x_k)$ and $u_k'=u_k(x_k+r_k x)$, we get 
$$
\rho(u_k',x,D_R,\lambda)>\frac{1}{2},\s \forall x\in D_{\frac{R}{2}}.
$$

Then  $\nabla u_k$ is bounded in $L^{2,\infty}(D_R)$ for any $R$. Set $c_k$ to be mean value of $u_k$ over $D$. By Poincare inequality,
${u_k}'-c_k$ is bounded in $W^{1,p}(D_R)$. Then we may assume 
${u_k}'-c_k$ converges to $u$ weakly in $W^{1,p}_{loc}$,
where
$$\left\{\begin{array}{l}
-\Delta u=0,\\ [\mv]
\frac{1}{|D|}\int_{|D|}|u-\frac{1}{|D|}\int_D u|dx=\lambda,\\[\mv]
 \frac{1}{|D(x)|}\int_{D(x)}|u-\frac{1}{|D(x)|}\int_{D(x)}u|dx\leq \lambda.
\end{array}\right.
$$ 
Hence, we get
$$
\|u-u(x_0)\|_{L^1(D(x_0))}\leq\lambda,\s \forall x_0,
$$
which implies that $\|\nabla u\|_{L^\infty(D_\frac{1}{2}(x_0))}<C$. Then
$\nabla u$ is a constant vector. 

Since $\int_D|u-u(0)|=\lambda|D|$, $u$ can not be  a constant.
\endproof

Theorem \ref{main1} is a corollary of the following lemma.
\begin{lem}\label{vol1}
Let $g=
e^{2u}g_{euc}$ be a metric defined on $D$, and $K_{g}$ be the Gauss
curvature. We assume there exists  $r_0>0$ and $\Lambda_1$, such that 
\begin{equation}\label{density}
\frac{\vol(B_r^g(x)\cap D_\frac{3}{4},g)}{\pi r^2}<\Lambda_1
\end{equation}
for any $x\in D_\frac{1}{2}$ and $r<r_0$. Then, there exists 
$\epsilon_0=\epsilon_0(\Lambda_1)$ and $a=a(\epsilon_0,r_0,\Lambda_1)$, such that if $\int_D|K_{g}|d\mu_g<\epsilon_0$, then
$$
\inf_{D_\frac{1}{2}}\rho(u,x,D,\lambda)>a>0.
$$
\end{lem}

\proof 
Assume this is not true.  By Lemma \ref{linear.function} and Corollary \ref{BH2}, we can find
$x_k\rightarrow x_0\in \overline{D_\frac{1}{2}}$ and $c_k$ such that $u_k(x_k+r_kx)-c_k
\rightarrow u_0=ax^1+bx^2+c$ weakly in $W^{1,p}_{loc}(\R^2)$, and
$e^{u_k(x_k+r_kx)-c_k}$ converges in $L^q_{loc}(\R^2)$ for any $q$. For simplicity, we assume $u_0=x^1$ and define $g_0=e^{2u_0}
(dx^1\otimes dx^1+dx^2\otimes dx^2)$.

Set $T(\theta)$ to be the constant, such that 
$$
Length((\cos\theta,\sin\theta)t|_{t\in[0,T(\theta)]},g_0)=R.
$$
It is easy to check that
$$
T(\theta)=\frac{\log (1+R\cos\theta)}{\cos\theta},\s i.e.\s
e^{T(\theta)\cos\theta}=R\cos\theta+1
$$
Let $a>0$ be sufficiently small and put
$$
\Omega(R)=\{(r,\theta):r\in(0,T(\theta)),\s \theta\in(-\frac{\pi}{2}+a,\frac{\pi}{2}-a)\}.
$$
Since
\begin{eqnarray*}
\vol(\Omega(R),g_0)
&=&
\int_{-\frac{\pi}{2}+a}^{\frac{\pi}{2}-a}\left.\left(
\frac{e^{2r\cos\theta}r}{2\cos\theta}-\frac{e^{2r\cos\theta}}{4\cos^2\theta}\right)\right|_0^{T(\theta)}d\theta\\
&=&
\int_{-\frac{\pi}{2}+a}^{\frac{\pi}{2}-a}\left(\frac{(R\cos\theta+1)^2T(\theta)}{2\cos\theta}-\frac{(R\cos\theta+1)^2-1}{4\cos^2\theta}\right)d\theta\\
&=&R^2\int_{-\frac{\pi}{2}+a}^{\frac{\pi}{2}-a}\left(\frac{1}{2}\log(1+R\cos\theta)
-\frac{1}{4}\right)d\theta
+\int_{-\frac{\pi}{2}+a}^{\frac{\pi}{2}-a}RT(\theta)d\theta\\
&&+\int_{-\frac{\pi}{2}+a}^{\frac{\pi}{2}-a}\frac{T(\theta)-R}{2\cos\theta}d\theta\\
&=&R^2\int_{-\frac{\pi}{2}+a}^{\frac{\pi}{2}-a}\left(\frac{1}{2}\log(1+R\cos\theta)
-\frac{1}{4}\right)d\theta
+O(R\log R),
\end{eqnarray*}
we have
\begin{equation*}
\lim_{R\rightarrow+\infty}\frac{\vol(\Omega(R),g_0)}{\pi R^2}=+\infty.
\end{equation*}

Let $T_1=\max_{S^1} T(\theta)$.  We have
$$
\int_{0}^{2\pi}\int_{0}^{T_1}|e^{u_k(r,\theta)}-e^{u_0(r,\theta)}|drd\theta=
\int_{0}^{2\pi}\int_{0}^{T_1}|e^{u_k(r,\theta)}-e^{u_0(r,\theta)}|r^\frac{1}{p}r^{-\frac{1}{p}}drd\theta\leq
C\|e^{u_k}-e^{u_0}\|_{L^p}\rightarrow 0.
$$
After passing to a subsequence, we can find $A\subset S^1$, such that
$L^1_{S^1}(A)<\epsilon$ and  
$$
\int_0^{T(\theta)}|e^{u_k}-e^{u_0}|dr\rightarrow 0, \forall \theta\notin A.
$$
Set
$$
\Omega(R,A)=\Omega(R)\setminus \{(r,\theta):\theta\in A\}.
$$
Then, we have 
$$
\Omega(R,A)\subset B_{R+1}^{g_k}(0)
$$
when $k$ is sufficiently large. Hence, 
we can choose $\epsilon$ small enough such that 
$$
\vol(\Omega(R,A),g_k)\geq
\frac{1}{2}\vol(\Omega(R),g_0).
$$
Then for any $K>0$, we can find $R$, such that 
$$
\frac{\vol(B_{R+1}^{g_k}(0),g_k)}{\pi(R+1)^2}>K
$$
when $k$ is sufficiently large. Then we get
$$
\frac{\vol(B_{(R+1)r_k}^{g_k}(x_k),{g_k})}{\pi((R+1)r_k)^2}>K.
$$
This is a contradiction.
\endproof

{\it The proof of Corollary \ref{main2}: } It is well-known that
on a   complete Riemann surface$(\Sigma,g)$, 
\begin{equation}\label{volume.comparison}
\frac{\vol(B^g_r(p),g)}{\pi r^2}<1+\int_{B^g_r(p)}K^{-}.
\end{equation}
In fact, if we let $\Omega\subset
T_x\Sigma$ be the segment domain, and 
$$
exp_x^*(g)=dr\otimes dr+\Theta^2 d\theta\otimes d\theta,
$$
then by the Jacobi equation $\Theta_{rr}+K\Theta=0$ in $\Omega_p$, 
we have
$$
\Theta=r-\int_0^r\int_0^tK\Theta(\tau,\theta)d\tau dt
\leq r+\int_0^r\int_0^tK^{-}\Theta d\tau dt.
$$
Let 
$$
\hat{\Theta}=\left\{\begin{array}{ll}
\Theta&x\in\Omega\\
0&x\notin\Omega.
\end{array}\right.
$$
We get
$$
\vol(B^g_r(x),g)=\int_{B^{exp^*_x(g)}_r\cap exp^{-1}_x(D)}\hat{\Theta}\leq
\int_0^{r}\int_0^{2\pi}\hat{\Theta}drd\theta\leq \pi r^2+r^2\int_{B^g_r}K^-d\mu_g.
$$

\endproof

We will use the following simple lemma.
\begin{lem}\label{extension}
If $u$ is smooth and $\|u\|_{W^{1,1}(D)}+\|\Delta u\|_{L^1(D)}<\Lambda_3$, then $e^{2u}g_{euc}|_{D_\frac{1}{2}}$ can be extended to a complete metric on $
\C$ with $\|K\|_{L^1}<C(\Lambda_3)$. Moreover, we have
$$
\|\nabla u\|_{L^{2,\infty}(D_\frac{1}{2})}<C(\Lambda_3).
$$
\end{lem}

\proof
Let $\eta$ be a cut-off function which is 1 on $D_\frac{5}{8}$ and 0 on $D^c_\frac{7}{8}$. Set $g'=e^{2\eta u}g_{euc}$. We have
$$
-\Delta \eta u=-u\Delta\eta-2\nabla u\nabla\eta-\eta \Delta u.
$$
Then 
$$
K(g')=e^{-2\eta u}(-u\Delta\eta-2\nabla u\nabla\eta-\eta \Delta u).
$$
Hence
$$
\int_\C|K(g')|d\mu_{g'}\leq C\int_D(|u|+|\nabla u|+|\Delta u|)<C(\Lambda_3).
$$
\endproof

\section{convergence of distance function}

In this section, we  set $g_k=e^{2u_k}g_{euc}$ which satisfies
\begin{itemize}
\item[{\rm 1)}]
$u_k=0$ on $D_2^c$;
\item[{\rm 2)}]
$u_k$ converges weakly to $u_0$ in $W^{1,p}(D_2)$ for any $p\in (1,2)$;
\item[{\rm 3)}]
$
\int_D|K(g_k)|d\mu_{g_k}<\epsilon_0.
$
\end{itemize}
By the results in the last section, we have 
$$
\|\nabla u_k\|_{L^{2,\infty}(D)}<C.
$$
Without loss of generality, we assume $|K(g_k)|d\mu_{g_k}$ and $K^{-}(g_k)d\mu_{g_k}$ converge to a measure $\nu$ and $\nu^-$
respectively in the sense of distribution. 

Let $d_k$ be the distance function defined by $g_k$. Then we have
$$
|e^{-u_k(x)}\nabla_xd_k(x,y)|=|e^{-u_k(y)}\nabla_yd_k(x,y)|=1.
$$
Therefore, 
$$
\int_{D\times D}(|\nabla_x d_k|^q+|\nabla_y d_k|^q)dxdy<C,
$$
where $q<\frac{4\pi-\epsilon_0}{\epsilon_0}$.
Hence, $d_k(x,y)$ converges to a function $d_0$ in $C^0(D
\times D)$. 

First, we prove the following lemma.
\begin{lem}\label{GH.smooth}
If $u_0=0$ and $\nu^-=0$,
 then 
$$
d_0(0,x)=d_{g_{euc}}(0,x),\s \forall x\in D.
$$
\end{lem}

\proof
Note that $\nabla e^{2u_k}=2e^{2u_k}\nabla u_k$, which is bounded
in $L^q$ for some $q>1$.
By the trace embedding theorem, 
it is easy to check that $d_0\leq d_{g_{euc}}$, therefore, we only need to show
$d_0\geq d_{g_{euc}}$.  Assume there exists $x'$, such that
$$
r=d_0(0,x')< d_{g_{euc}}(0,x').
$$
Then $x'\notin \overline{D_r}$,   $B^{d_0}_r(0)\setminus \overline{D_r}$ ,which  is a non-empty open set, hence
$$
\vol(B^{d_0}_r(0)\setminus \overline{D_r},g_{euc})>a>0.
$$
Here $B^{d_0}_r(0)=\{x:d_0(0,x)<r\}$.
 Since $\|K^{-}(g_k)\|_{L^1}\rightarrow 0$, we get
$$
\varlimsup_{k\rightarrow+\infty}\frac{\vol(B_t^{g_k}(0),{g_k})}{\pi t^2}\leq 1,\s \forall t>0.
$$

Since
$B_{r-\epsilon}^{d_0}(0)\subset B_{r}^{g_k}(0)$
when $k$ is sufficiently large, we have
$$
\vol(B_{r-\epsilon}^{d_0}(0),{g_{euc}})=
\lim_{k\rightarrow+\infty}\vol(B^{d_0}_{r-\epsilon}(0),{g_k})\leq \varlimsup_{k\rightarrow+\infty}\vol(B^{g_k}_r(0),{g_k})\leq
\pi r^2.
$$
Let $\epsilon\rightarrow 0$, we get a contradiction.
\endproof

As an application, we  give the proof of Theorem \ref{neck.estimate}.

{\it The proof of Theorem \ref{neck.estimate}:}
Assume the left inequality of \eqref{diam} is no true.  Then we can find $u_k$ with $\|\nabla u_k\|_{L^{2,\infty}(D)}<\Lambda_3$, and
$\theta_k$, $\theta_k'\in[0,2\pi)$, such that
\begin{equation*}
\frac{d_{g_k}(e^{i\theta_k},2e^{i\theta'_k})}{e^{c_k}}\rightarrow 0. 
\end{equation*}
Let $c_k$ be  the mean value of $u_k$ on $\partial D_\frac{3}{2}$.
By Theorem \ref{main1} and Corollary \ref{BH2},  $u_k-c_k$ converges to a harmonic function with $\int_{\partial D_\frac{3}{2}}w=0$, and $\|\nabla w\|_{L^p}<C(\Lambda_3)$. Then we get
$$
\|w\|_{C^0}<C(\Lambda_3).
$$
Let $g_k'=e^{2u_k-2c_k-2w}g_{euc}$. We have
$$
\lim_{k\rightarrow+\infty}d_{g_k'}(e^{i\theta_k},2e^{i\theta'_k})=|e^{i\theta}-2e^{i\theta'}|,
$$
where $\theta$ and $\theta'$ are the limits of $\theta_k$ and 
$\theta_k'$ respectively.
Note that 
$$
\frac{d_{g_k}(e^{i\theta_k},2e^{i\theta'_k})}{e^{c_k}}\geq C(\Lambda_3) d_{g_k'}(e^{i\theta_k},2e^{i\theta'_k}).
$$
We get a contradiction. 

The proofs of the other parts of this theorem are similar, hence are omitted here. \endproof

%Next, we prove the following

\begin{lem}\label{distance}
$d_0$ is a distance function. Moreover, for any $\delta>0$,
we can find $a(\delta)>0$, such that 
$$
d_0(x,y)>a(\delta),\s whenever\s |x-y|\geq\delta.
$$
Or equivalently, $\phi(x,y)=|x-y|$ is continuous on $(\C,d_0)$ 
\end{lem}

\proof Without loss of generality,   assume $|x-y|=\delta$.
We may choose $m$, such that $\frac{\int_\C|K(g_k)|d\mu_{g_k}}{m}<\epsilon_2$. Then, after passing to a subsequence, we can choose
$i<m$, such that 
$$
\int_{D_{2^{-i}\delta}\setminus D_{2^{-i-1}\delta}(x)}|K(g_k)|e^{2u_k}
<\epsilon_2
$$ 
for any $k$. By the trace embedding theorem, we may assume that 
$$
|\int_{\partial D_{2^{-i}\delta}(x)}u_k|<C(\delta).
$$
By Theorem \ref{neck.estimate}, we have 
$$
d_{g_k}(x,y)\geq\lambda(\delta)>0.
$$
Thus, $d_0(x,y)=0$ implies that $x=y$.
\endproof

Now, we start to proof $d_0=d_{u_0}$.
By the trace embedding theorem, we have
$$
\int_{\gamma}e^{u_k}\rightarrow\int_\gamma e^u.
$$
Then, we get
$$
d_0(x,y)\leq d_{u_0}(x,y).
$$
Thus, it only needs to check whether  $d_0(x,y)\geq d_{u_0}(x,y)$.
The key observation is the following result.
\begin{lem}\label{du0/d0}
For any $\epsilon$, we can find $\beta$ and $\tau$, such that if 
$$
\nu(B_{\delta}(x))<\tau, \s \delta<\frac{1}{2},\s x\in D_\frac{1}{2},
$$
then
$$
\frac{d_{u_0}(x,y)}{d_0(x,y)}\leq1+\epsilon,\s \forall y\in D_{\beta\delta}(x).
$$
\end{lem}

\proof 
Assume the result is not true. Then we can find $\delta_m\in (0,\frac{1}{2})$,  $x_m\in D_\frac{1}{2}$, $y_m\in D$, such that
$\frac{|x_m-y_m|}{\delta_m}\rightarrow 0$ and  
$$
\lim_{m\rightarrow+\infty}
\limsup_{k\rightarrow+\infty}\int_{D_{\delta_m(y_m)}}|K(g_k)|d\mu_{g_k}= 0,\s \frac{d_0(y_m,x_m)}{d_{u_0}(y_m,x_m)}\rightarrow l_0\leq\frac{1}{1+\epsilon}.
$$
For any fixed $m$, we can find $k_m$, such that 
$$
\left|\frac{d_{g_{k_m}}(y_m,x_m)}{d_{u_0}(y_m,x_m)}-\frac{d_{0}(y_m,x_m)}{d_{u_0}(y_m,x_{m})}\right|<\frac{1}{m},
$$
and
\begin{equation}\label{L1}
\frac{1}{|D_{r_m}|}\int_{D_{r_m}(y_m)}|u_{k_m}-u_0|<\frac{1}{m},
\end{equation}
where $r_m=|y_m-x_m|$.
For simplicity, we set $x_m=y_m+r_m(0,1)$ and $u_m'=u_{k_m}(y_m+r_mx)-c_k$, where $c_k$ is chosen such that
$$
\int_{D_1}u_m'=0.
$$
We set $g_m'=e^{2u_m'}g_{euc}$.  By Corollary \ref{main2},  $\|\nabla u_m'\|_{L^{2,\infty}(D(z)}<C$ for any $z$, then we may assume $u_m'$ converges to a harmonic function $u$ weakly in $W^{1,p}_{loc}(\C)$.
Moreover, for any $D_r(z)$, we have 

\begin{eqnarray*}
r^{2-q}\int_{D_r(z)}|\nabla u|^q&\leq&\lim_{k\rightarrow+\infty}
(rr_m)^{2-q}\int_{D_{rr_m}(y_m+r_mz)}|\nabla u_{k_m}|^q\\
&\leq& C\|\nabla u_{k_m}\|^q_{L^{2,\infty}(D_{rr_m}(y_m+r_mz))}\\
&<&\Lambda.
\end{eqnarray*}
Then $u$ is a constant with $\int_{D}u=0$. Hence,
$u=0$.
By \eqref{L1},  $(u_0(x_m+r_mx)-u_{k_m}(x_m+r_mx))$ converges to 0
in $W^{1,p}_{loc}(\C)$weakly, and
$$
e^{-c_k}d_{u_0}(x_m,y_m)\leq \int_{[0,1]}e^{u_0(y_m+r_mx)-c_k}=\int_{[0,1]}e^{(u_0(y_m+r_mx)-u_{k_m}(y_m+r_mx))
+u_m'}.
$$
By the trace embedding theorem, we get
$$
\lim_{k\rightarrow+\infty} e^{-c_k}d_{u_0}(x,y)\leq \int_{[0,1]}e^{u_0-c_k}\leq 1
$$
By Lemma \ref{GH.smooth},
$$
e^{-c_k}d_{g_{k_m}}(x_m,y_m)=d_{g_m'}(0,(0,1))\rightarrow 1.
$$
Then
$$
l_0>\frac{d_0(x,y)}{|d_{u_0}(x,y)|}\geq1,
$$
which is impossible.
\endproof

\begin{pro}
$d_0=d_{u_0}$.
\end{pro}

\proof Let $\epsilon$, $\tau$ and $\beta$  be as in Lemma \ref{du0/d0}.
We set $A_\tau=\{x:\nu(\{x\})>\tau\}$. Obviously, $A_\tau$
is a finite set.   Then, for any $\delta>0$ and $B_{2\delta}(x)
\cap A_\tau=\emptyset$,  we have
$$
\int_{B_\delta(x)}|K(g_k)|d\mu_{g_k}<\tau,
$$
when $k$ is sufficiently large. Then
$$
\frac{d_{u_0}(x,y)}{d_0(x,y)}<1+\epsilon
$$
whenever $|x-y|<\beta\delta$ and $x\notin B_\delta(A_\tau)$.

Let $\gamma$ be the segment defined in $(\C,d_0)$ connecting
$x_1$ and $x_2$, i.e. $\gamma:[0,a]\rightarrow (\C,d_0)$
is a continuous map which satisfies
$$
d_0(\gamma(s),\gamma(s'))=|s-s'|,\s\forall s,s'\in[0,a].
$$ 

First, we consider  the case when $\gamma\cap A_\tau=\emptyset$. We may assume
$$
d_{g_{euc}}(A_\tau,\gamma[0,a])>\delta>0.
$$
By Lemma \ref{distance}, we can find 
$$
s_0=0<s_1<\cdots<s_m=a,
$$
such that 
$$
|\gamma(s_{i+1})-\gamma(s_i)|<\beta\delta.
$$
Then
\begin{eqnarray}\label{d0.d0'}\nonumber
d_0(x_1,x_2)&=&\sum_{i=0}^m d_0(\gamma(s_i),\gamma(s_{i+1}))\\ 
&\geq& (1+\epsilon)^{-1}\sum_{i=0}^{m-1}d_{u_0}(\gamma(s_i),
\gamma(s_{i+1}))\\\nonumber
&\geq& (1+\epsilon)^{-1}d_{u_0}(x_1,x_2).
\end{eqnarray}

Now,  we consider the case in which $\gamma\cap A_\tau\neq
\emptyset$.  Let 
$$
\gamma\cap A_\tau=\{\gamma(a_1), \cdots, \gamma(a_i)\}.
$$ 
Then we have
\begin{eqnarray*}
d_0(x_1,x_2)&\geq& d_0(x_1,\gamma(a_1-\epsilon'))
+d_0(\gamma(a_1+\epsilon'),\gamma(a_2-\epsilon'))+
\cdots+d_0(\gamma(a_i+\epsilon',x_2))\\
&\geq&
(1-\epsilon)^{-1}(d_{u_0}(x_1,\gamma(a_1-\epsilon'))
+
\cdots+d_{u_0}(\gamma(a_i+\epsilon'),x_2)).
\end{eqnarray*}
Let $\epsilon'\rightarrow 0$, we get \eqref{d0.d0'} again.

Now, let $\epsilon\rightarrow 0$, we get the desired result.
\endproof

\section{Some applications}
In this section, we give some applications.

\subsection{A sequence defined on $D_r\setminus D_{\frac{r_k}{r}}$}
Let $g_k=e^{2u_k}g_{euc}$ be a metric defined on $D_r\setminus D_{\frac{r_k}{r}}$, where $r_k\rightarrow 0$. We assume 
\begin{itemize}
\item[{\rm 1)}] $(D_r\setminus D_{\frac{r_k}{r}},g_k)$ can be
extended to a complete surface with 
$$\|K(g_k)\|_{L^1}<C;
$$
\item[{\rm 2)}] 
$
\int_{D_{2t}\setminus D_t}|K(g_k)|e^{2u_k}<\epsilon_2,\s \forall t\in(\frac{r_k}{r},r);
$
\item[{\rm 3)}] 
$
\lim\limits_{r\rightarrow 0}\lim\limits_{k\rightarrow+\infty}\sup_{t\in[\frac{r_k}{r},r]}d_{g_k}(\partial D_{2t},\partial D_t)=0.
$
\end{itemize}

By Theorem \ref{neck.estimate}, 
$$
\frac{\vol(D_{2t}\setminus D_t,g_k)}{d^2_{g_k}(\partial D_{2t},\partial D_t)}<C,\s
\forall t\in[\frac{r_k}{r},r].
$$
Without loss of generality, we set $r=2^{m}\frac{r_k}{r}$. Then
\begin{eqnarray*}
\vol(D_r\setminus D_{\frac{r_k}{r}},{g_k})&=&\sum_{i=1}^m \vol(D_{2^{i}\frac{r_k}{r}}\setminus
D_{2^{i-1}\frac{r_k}{r}},{g_k})\\
&\leq& \sum_{i=1}^md^2_{g_k}(\partial D_{2^{i}\frac{r_k}{r}},\partial D_{2^{i-1}\frac{r_k}{r}})\\
&\leq& \epsilon \sum_{i=1}^md_{g_k}(\partial D_{2^{i}\frac{r_k}{r}},\partial D_{2^{i-1}\frac{r_k}{r}}).
\end{eqnarray*}
Thus, when $diam_{g_k}(D_{r}\setminus D_{\frac{r_k}{r}})$ is bounded  above, we get
\begin{equation}\label{neck.volume}
\lim_{r\rightarrow 0}\lim_{k\rightarrow+\infty}\vol(D_r\setminus D_{\frac{r_k}{r}},{g_k})=0.
\end{equation}
As an application, we have the following result.
\begin{pro}\label{bubble}
Assume
$$
\int_{D}(1+|K(g_k)|)e^{2u_k}<C,\s \|\nabla u_k\|_{L^{2,\infty}(D)}<C,
$$
and
$$ 
\lim_{\delta\rightarrow 0}\lim_{k\rightarrow+\infty}\vol(D_\delta,g_k)>0.
$$
Then we can find $x_k\rightarrow x_0$, $r_k\rightarrow 0$, and a finite set $\mathcal{S}$,
such that 
$u_k(x_k+r_kx)-\log r_k$ converges  weakly in $W^{1,p}_{loc}(\C\setminus \mathcal{S})$ to a function $u$ with
$$
\int_{\C}e^{2u}<+\infty, \s where\s p\in(1,2).
$$
\end{pro}

\proof Let $\epsilon_0'=\min\{\epsilon_0,\epsilon_2\}$.
We assume $\int_D|K(g_k)|e^{2u_k}\leq m\frac{\epsilon_0'}{2}$. We will prove the result by induction
of $m$.
When $m=1$, 
we set 
$$
r_k(x)=\sup\{t:\int_{D_t(x)}|K(g_k)|e^{2u_k}\leq\frac{\epsilon_0'}{2}\},
$$
and take $x_k$ such that $\int_{D_{r_k}(x_k)}|K(g_k)|e^{2u_k}=\frac{\epsilon_0'}{2}$. If there exists 
$r_k'\rightarrow 0$, $\frac{r_k'}{r_k}\rightarrow+\infty$, such that
$$
d_{g_k}(\partial D_{2r_k'}(x_k),\partial D_{r_k'}(x_k))\rightarrow \lambda>0,
$$
then by Theorem  \ref{neck.estimate}, $\vol(D_{2r_k'}(x_k)\setminus
 D_{r_k'}(x_k),g_k)>\lambda'>0$. By Theorem \ref{main1},
$u_k(r_k'x+x_k)$ converges weakly in $W^{1,p}(\C\setminus\{0\})$. 

Now, we assume 
$$
\lim_{r\rightarrow 0}\lim_{k\rightarrow+\infty}\sup_{t\in[\frac{r_k}{r},r]}d_{g_k}(\partial D_{2t}(x_k),\partial D_t(x_k))=0.
$$
We claim that $u_k'(x)=u_k(r_kx+x_k)$ must converge. Since $\|\nabla u_k'\|_{L^{2,\infty}}(D_R)<C(R)$, $u_k'-c_k$ must converge
weakly in $W^{1,p}(D_R)$, where $c_k$ is the mean value of $u_k'$
on $D$. Since $\int_De^{2u_k'}<C$, it follows from Jensen's inequality that  $c_k<+\infty$. Thus,
the fact that $u_k'$ does not converge implies that $c_k\rightarrow-\infty$,
hence $\vol(D_{\frac{r_k}{r}}(x_k))\rightarrow 0$, which implies 
$$
\lim_{r\rightarrow 0}\lim_{k\rightarrow+\infty}\vol(D_r(x_0),g_k)
=\lim_{r\rightarrow 0}\lim_{k\rightarrow+\infty}\vol(D_r(x_k)\setminus D_{\frac{r_k}{r}},g_k)=0.
$$
We get a contradiction.

Now, we
assume $\int_D|K(g_k)|e^{2u_k}\leq (m-1)\frac{\epsilon_0'}{2}$ implies the result and
set $\int_D|K(g_k)|e^{2u_k}\leq m\frac{\epsilon_0'}{2}$. We assume $u_k'$ does not converge and
set 
$$
t_k=\sup\{t:\int_{D_t(x)}|K(g_k)|e^{2u_k}\leq\frac{\epsilon_0'}{2}(m-1)\}.
$$
Assume $\int_{D_{t_k}(x_k)}|K(g_k)|e^{2u_k}=\frac{\epsilon_0'}{2}(m-1)$. If there exists 
$t_k'\rightarrow 0$, $\frac{t_k'}{t_k}\rightarrow+\infty$, such that
$$
d(\partial D_{2t_k'}(x_k),\partial D_{t_k'}(x_k))\rightarrow \lambda>0,
$$
then $u_k(t_k'x+x_k')-\log t_k'$ converges weakly in $W^{1,p}(\C\setminus\{0\})$.

Now, we assume 
$$
\lim_{r\rightarrow 0}\lim_{k\rightarrow+\infty}\sup_{t\in[\frac{t_k}{r},r]}d(\partial D_{2t}(x_k),\partial D_t(x_k))=0,
$$
which implies that 
$$
\lim_{r\rightarrow 0}\lim_{k\rightarrow+\infty}\vol(D_{r}\setminus D_\frac{t_k}{r}(x_k))=0.
$$
Since $u_k'$ does not converge, we have $\frac{t_k}{r_k}\rightarrow 0$. Otherwise, we get
$$
\lim_{r\rightarrow 0}\lim_{k\rightarrow+\infty}\vol(D_{r}(x_k),g_k)=0.
$$
Put $u_k''=u_k(x_k+t_kx)-\log t_k$. Then, we have
$$
\lim_{k\rightarrow+\infty}\vol(u_k'',D_R)=\lim_{k\rightarrow+\infty}\vol(u_k,D_{Rt_k}(x_k))>0.
$$
Without loss of generality, we assume 
$K(x_k+t_kx)e^{2u_k''}dxdy$ converges to $\nu$ in the sense of
distribution. 
Choose $c_k$ such that  $u_k''-c_k$ converges weakly in 
$W^{1,p}(\C\setminus \mathcal{S}'')$, where 
$$
\mathcal{S}''=\{z:\nu(\{z\})\geq\frac{\epsilon_0'}{2}\}.
$$
If $c_k$ is bounded, then $u_k''$ converges weakly. Thus we may assume $c_k
\rightarrow-\infty$, which implies
$$
\sum_{z\in\mathcal{S}''}\lim_{k\rightarrow+\infty}
\vol(D_r(z),g_k)>0.
$$ 
However, $r_k/t_k\rightarrow 0$ implies that $\mu(\{0\})\geq
\frac{\epsilon_0'}{2}$. Together with $\nu(D_1^c)\geq
\frac{\epsilon}{2}$, we have
$$
\mu(\{z\})\leq\frac{m-1}{2}\epsilon_0',\s\forall z\in\mathcal{S}''.
$$
Using the induction on $u_k''$, we will get the result.
\endproof

\subsection{A Mumford type lemma}
Let $(\Sigma,g)$ be a Riemann surface and $\sigma\in\pi_1(\Sigma)$. We denote the length of the shortest closed geodesic representing $\sigma$ by $L^-(\sigma)$. We will use the following  lemmas (see \cite{DeTurck-Kazdan} and
\cite{Hummel} for proofs):
\begin{lem}\label{isothermal}
Let $g_k,g$ be smooth Riemannian metrics on a surface $M$,
such that $g_k \to g$ in $C^{s,\alpha}(M)$, where $s \in N$,
$\alpha \in (0,1)$. Then for each $p \in M$ there exist
neighborhoods $U_k, U$ and smooth conformal 
diffeomorphisms
$\varphi_k:D \to U_k$, such that $\varphi_k \to \varphi$
in $C^{s+1,\alpha}(\overline{D},M)$.
\end{lem}
\begin{lem}\label{conformal.class}
If $Conf(g_k)$ converges, then there exists $h_k$ which is conformal to
$g_k$, such that $h_k$ converges smoothly. 
\end{lem}

First, we prove the following:
\begin{thm}\label{Mumford}
Let $(\Sigma,g_k)$ be a Riemann surface  with genus $\geq1$.
We assume $\|K(g_k)\|_{L^1}=1$. Then if the conformal class $Conf(g_k)$ is induced by $g_k$ diverges, then there exists
a nontrivial $\sigma\in\pi_1(\Sigma)$ such that $L^-(\sigma)
\rightarrow 0$. Conversely, if there exist $\sigma_1$ , $\sigma_2\in\pi_1(\Sigma)$ and some number $l_0$, such that  $L^-(\sigma_1,g_k)\rightarrow 0$ and $L^-(\sigma_2,g_k)>l_0>0$, then $Conf(g_k)$ diverges. 
\end{thm}

\proof
First of all, we assume $Conf(g_k)$ diverges. By the Collar Lemma, we can
find $T_k\rightarrow+\infty$ and $\Omega_k\subset\Sigma$, such that $\Omega_k$
is conformal to $S^1\times[-T_k,T_k]$, where the homotopy class of $S^1\times\{0\}$ is nontrivial in $\pi_1(\Sigma)$. Obviously, we can find
$a_k<\frac{T_k}{2}$, such that 
$$
\int_{S^1\times[a_k-1,a_k+1]}(1+|K(g_k)|)e^{2u_k}\rightarrow 0.
$$
Put $g_k=e^{2u_k}(dt^2+d\theta^2)$. By Corollary \ref{collapse},
$l_{g_k}(S^1\times\{a_k\})\rightarrow 0$.

Next, we will show  when $Conf(g_k)$ converges, if
there exists nontrivial $\sigma\in\pi_1(\Sigma)$, such that $L^-(\sigma,g_k)\rightarrow 0$, then $L^-(\sigma',g_k)
\rightarrow 0$ for any nontrivial $\sigma'$.
 
First, we consider the case when  $Conf(g_k)$ is fixed.
Assume $g_k=e^{2u_k}g$ for a smooth metric $g$.
Without loss of generality,  assume $\nu$ be limit measure of $|K(g_k)|e^{2u_k}dx$ in the sense of
distribution and set 
$$
\mathcal{S}(\{g_k\})=\{x:\nu(\{x\})>\frac{\epsilon_0}{2}\}.
$$
Set $\gamma_k:[0,a_k]\rightarrow \Sigma$ to be the shortest nontrivial closed geodesic with the  unit speed, which represents $\sigma'$. For simplicity, we assume 
$\gamma_k(0)\rightarrow x_0$. Select an isothermal coordinate
system $(D,(x^1,x^2))$ around $x_0$. Let $2m\epsilon_0<1$. Since $\sigma'$
is nontrivial in $\pi_1$, we can find $b_0=0<b_1\leq b_2<\cdots<b_{2m+1}=a_k$,
such that $\gamma_k([b_{2i},b_{2i+1}])\subset\overline{ D_{(i+1)/(2m)}\setminus
D_{i/(2m)}}$ and $\gamma_k(b_{2i+1})$, $\gamma_k(b_{2i+2})\subset
\partial D_{(i+1)/(2m)}$, where $i=1,\cdots,m$.  Choose $i_0$, such that
$$
\int_{D_{(i_0+1)/(2m)}\setminus D_{(i_0-2)/(2m)}}|K(g_k)|e^{2u_k}<\epsilon_0,
$$
and
let $c_k$ be the mean value of $u_k$ on a small disk contained in $D_{i_0}\setminus D_{i_0-1}$, which contains no
Condensation points. By Theorem \ref{main1},  $\|\nabla u_k\|_{L^{2,\infty}}$ is bounded on any $\Omega
\subset\subset\Sigma\setminus \mathcal{S}(\{g_k\})$. Then, by 
Poincar\'{e} inequality, $\|u_k-c_k\|_{W^{1,p}(\Omega)}$ is bounded when 
$$
D_{i_0/(2m)}\setminus D_{(i_0-1)/(2m)}\subset\Omega\subset\subset\Sigma\setminus \mathcal{S}(\{g_k\}).
$$
Thus we must have $c_k\rightarrow -\infty$,
for $d(D_{i_0},D_{i_0-1},g_k)\leq l_{g_k}(\gamma_k)\rightarrow 0$. 
This implies that $e^{u_k}$ converges to 0  in $W^{1,p}_{loc}(\mathcal{S}(\{g_k\}))$.
 Then $L^-(\gamma)\rightarrow 0$
for any $\gamma\neq 1$ in $\pi_1(\Sigma)$. We get a contradiction.

Now, we assume $Conf(g_k)$ converges. By Lemma \ref{isothermal}
and \ref{conformal.class}, the proof for this case is almost the same with the case when $Conf(g_k)$ is fixed,  hence is omitted here.
\endproof

\begin{cor}
Let  $g_k$ and $g_\infty$ be  metrics on a  closed Riemann surface. We assume $\|K(g_k)\|_{L^1}<C$ and $(\Sigma,g_k)$
converges to $(\Sigma,g_\infty)$ in the sense of Gromov-Hausdorff distance. Then 
$\vol(\Sigma,g_k)\rightarrow \vol(\Sigma,g_\infty)$.
\end{cor}

\proof Note that it follows from \eqref{volume.comparison} that $\mu(g_k)<C$.

First, we prove $Conf(g_k)$ converges. Assume this is not true.
By the proof of the Theorem \ref{Mumford}, we can find $\Omega_k\subset \Sigma_k$, which is conformal to $(S^1\times [-T_k,T_k],
e^{2u_k}(dt^2+d\theta^2))$
with $T_k\rightarrow+\infty$, and $a_k<\frac{T_k}{2}$,
such that the homotopy class of $S^1\times\{0\}$ is nontrivial in $\pi_1$ and
$diam(S^1\times(a_k-1,a_k+1))\rightarrow 0$. In fact, by
Theorem \ref{main1}, \ref{main2}, and trace embedding
inequality, we may assume 
$\sup_{t\in [a_k-1,a_k+1]}l_{g_k}(S^1\times\{t\})\rightarrow 0$. 
We  replace
$g_k$ with $e^{2u_k}|a_k-t|(dt^2+d\theta^2)$ on $S^1\times[a_k-1,a_k+1]$, which can be considered as a metric defined
on the cone
$$
C_k=\{(|t|\cos\theta,|t|\sin\theta, t):t\in [a_k-1,a_k+1]\}
$$
with $diam(C_k)\rightarrow 0$. Denote the new metric space 
by $\Sigma_k'$.  By Theorem 2.1 in \cite{Sormani-Wei}, whose
proof can be also found in \cite[page 100]{Gromov} and \cite{Cassorla}, 
there exists a surjective homomorphism
from  $\pi_1(\Sigma_k')$ to $\pi(\Sigma)$. This is impossible.\\

By Lemma \ref{isothermal}
and \ref{conformal.class}, we may assume $Conf(g_k)$ is fixed.
Let $g_k=e^{2u_k}g$,
where $g$ is a fixed metric. 
We show that   $u_k$
converges weakly in $W^{1,p}(\Sigma\setminus \mathcal{S}(\{g_k\}))$.
Assume this is not true. Then there exists $c_k\rightarrow-\infty$, such that $u_k-c_k$  converges weakly in $W^{1,p}(\Sigma\setminus \mathcal{S}_0)$ for some finite set 
$\mathcal{S}_0$. Then we can find embedded curves $\gamma_1$,
$\cdots$, $\gamma_m$, such that 
$$
\gamma_i\cap\mathcal{S}_0=\emptyset,\s \lan \gamma_1,\cdots,\gamma_m\ran=\pi_1(\Sigma).
$$
By the trace embedding inequality,  $l_{g_k}(\gamma_i)\rightarrow 0$. Let $\varphi$ be a smooth function which is 0 when $t\leq 0$,
1 when $t>\delta_0$ and positive on $(0,+\infty)$, here 
$\delta_0<\frac{1}{2}d_g(\mathcal{S}_0,\cup_i\gamma_i)$. Set 
$d_0(x)=d(x,\cup_i\gamma_i)$ and $g_k'=g_k\varphi(d_0)$.
Then $g_k'$ define a distance function $d_k'$ on $\Sigma/\sim$, where we say
$x\sim y$ if $x,y\in \cup_i\gamma_i$. Obviously, we have
$$
d_{GH}((\Sigma/\sim,d_k'),(\Sigma,d_{g_\infty}))\rightarrow 0.
$$
By Theorem 2.1 in \cite{Sormani-Wei} again, $\pi_1(\Sigma)$ is trivial, which is impossible.\\

Next, we claim that
for any
$p\in \mathcal{S}(\{g_k\})$, it must hold true that
$$
\lim_{r\rightarrow 0}\lim_{k\rightarrow+\infty}\vol(B_r^{g_0}(p),g_k)= 0.
$$
Otherwise, by Lemma \ref{bubble}, we can find $x_k\rightarrow p$, $r_k\rightarrow 0$ such that $u_k'=u_k(x_k+r_kx)$ converges weakly in $W^{1,p}(\C\setminus \mathcal{S}')$, where $\mathcal{S}'$ is a finite set. 

Let $u'$ be the limit. Then we have
$\int_{D_{2R}\setminus D_R}e^{2u'}\rightarrow 0$ as $R\rightarrow+\infty$. Let $u_R'=u(Rx)-\log R $, and $g_{R}'=e^{2u_R'}g_{euc}$. Then $(D_2\setminus D_1,g_{R}')$ is a new parametrization of $(D_{2R}\setminus D_R,g)$.
By Corollary \ref{main2}, we have 
$$
\lim_{R\rightarrow+\infty}L(\partial D_{R},g_R)=\lim_{R\rightarrow+\infty}L(\partial D_{2},g_R')=0.
$$
Thus, we have
$$
\lim_{R\rightarrow+\infty}\lim_{k\rightarrow+\infty}
L(\partial D_{Rr_k}(x_k),g_k)=0.
$$
Then there must be a shortest closed geodesic $\gamma_k$ on $D_\frac{1}{R}(x_k)\setminus D_{Rr_k}(x_k)$, which is also nontrivial in $\pi_1(D_\frac{1}{R}(x_k)\setminus D_{Rr_k}(x_k))$, 
when $R$ and $k$
 are sufficiently large. In addition, the length of $\gamma_k$ converges to 0 as $k\rightarrow+\infty$ and $D_\frac{1}{R}(x_k)
 \setminus\{\gamma_k\}$ has just two connected components .
 
Without loss of generality, we assume $d_{GH}((\Sigma,g_k),(\Sigma,g_\infty))<2^{-k}$. We set $\Sigma_k=(\Sigma,g_k)$ and
$$
X=\left(\bigsqcup_{k=1}^\infty\Sigma_k\right)\bigsqcup\Sigma_\infty,
$$
 where $\bigsqcup$ is the disjoint union operator, and $d$ is an 
 admissible distance such that
$$
d_{H,X}(\Sigma_k,\Sigma_\infty)<2^{-k}.
$$
We set $x_\infty\in\Sigma_\infty$ to be the limit of a sequence $\hat{x}_k\in\gamma_k\subset X_k$.  Since $L(\gamma_k,g_k)\rightarrow 0$,
$\gamma_k$ converges to $x_\infty$.

Let $2r$ be the injective radius of $\Sigma_\infty$. Then $\overline{B_r(\hat{x}_k)}^{\Sigma_k}$ converges
to $\overline{B_r(x_\infty)}^{\Sigma_\infty}$ in the Hausdorff distance.
Since $l_{g_k}(\gamma_k)\rightarrow 0$, we have $\gamma_k\subset
\overline{B_r(\hat{x}_k)}^{\Sigma_k}$ for sufficiently large $k$. Then  $\overline{B_r(\hat{x}_k)}^{\Sigma_k}
\setminus\gamma_k$ has just two connected components. Let $y_k$,
$y_k'\in \partial B_{r-\epsilon}^{\Sigma_k}(\hat x_k)$ lie on different components,
which converge to $y_\infty$ and $y_\infty'$ respectively.
Since the
segment $\widehat{y_ky_k'}$ must pass through $\gamma_k$, the
segment $\widehat{y_\infty y_\infty'}$ must pass through $x_\infty$. In the same way, we can find $r'$, such that 
for any $y'\in B_{r'}(y_\infty')$
 in $\Sigma_\infty$, $\widehat{y_\infty y'}$  passes through $x_\infty$, which is impossible.
\endproof

\subsection{Metrics  On $S^2$ With   Small $\|K(g)-1\|_{L^1}$}
Let $g_{S^2}$ be the standard metric defined on $S^2$ with 
$K=1$ and $g_k=e^{2u_k}g_{S^2}$.
By Theorem A.1 in \cite{Nerz}, if $\|K(g_k)-1\|_{L^p}<\epsilon$ for some $\epsilon$ and $p>1$, then $(S^2,g_k)$ converges in $C^{1,\alpha}$. We will extend their result to the case   $p=1$. We first prove that
\begin{lem}\label{sphere}
Let $g_k=e^{2u_k}g_{\S2}$. Assume $\vol(g_k) \leq \Lambda$ and
$\|K(g_k)-1\|_{L^1(S^2,g_k)}\rightarrow 0$. After passing to a subsequence,
we can find a M$\ddot{o}$bius
transformation $\sigma_k$, such that $\sigma_k^*(g_k)$
converges to $g_{\S2}$ weakly in $W^{1,p}$. Moreover, we have
$\sigma_k^*(g_k)=e^{2u_k'}g_{\S2}$ with $u_k'$ converges to 0 in
$W^{1,p}$ for any $p\in(1,2)$ and $e^{u_k'}$ converges to 1
in $L^q$ for any $q>1$.
\end{lem}

\proof
Since 
\begin{eqnarray*}
\vol(g_k)&=&\int_{S^2}(1-K(g_k))d\mu_{g_k}+\int_{S^2}K(g_k)d\mu_{g_k}\\
&=& \int_{S^2}(1-K(g_k))d\mu_{g_k}+4\pi,
\end{eqnarray*}
we get 
\begin{equation}\label{area}
\vol(g_k)\rightarrow 4\pi.
\end{equation}
Then, we have
$$
\|K(g_k)\|_{L^1(S^2,g_k)}\leq\|K(g_k)-1\|+\vol(g_k)\rightarrow 4\pi.
$$

Assume  $|K(g_k)|d\mu_{g_k}$ converges to measure $\nu$ in the sense of
distribution. Put
$$
\mathcal{S}=\{x:\nu(\{x\})>\frac{\epsilon_0}{2}\}.
$$

Let $x_0\in\mathcal{S}$, $y_0$ be the antipodal point of $x_0$ on $\S2$ , and $\pi$ be the 
stereographic projection from $\S2\setminus\{y_0\}$ to $\C$. 
It is well-known that $\pi$
defines an isothermal coordinate system with $x_0=0$. 
In this new coordinate, we set 
$$
g_k=e^{2v_k}g_{euc}.
$$
We have  $-\Delta v_k=K_{{g}_k}e^{2v_k}$, and
$$
\lim_{r\rightarrow 0}\lim_{k\rightarrow+\infty}\vol(D_r,{g}_k)=\lim_{r\rightarrow 0}\lim_{k\rightarrow+\infty}\int_{D_r}|K_{{g}_k}|e^{2v_k}>\frac{\epsilon_0}{2}.
$$
We set 
$$
r_k(x)=\sup\{t:\int_{D_t(x)}|K({g}_k)|e^{2v_k}\leq\frac{\epsilon_0}{2}\},
$$
and take $x_k$ such that  $\int_{D_{r_k}(x_k)}|K({g}_k)|e^{2v_k}=\frac{\epsilon_0}{2}$. Put
$v_k'=v_k(x_k+r_kx)-\log r_k$, $g_k'=e^{2v_k'}g_{euc}$.
Since 
$$
\vol(D_R,g_k')=\int_{D_{Rr_k}(x_k)}e^{2v_k}\geq \int_{D_{Rr_k}(x_k)}|K(g_k)|e^{2v_k}-
\int_{D_{Rr_k}(x_k)}|K(g_k)-1|e^{2v_k}\rightarrow \frac{\epsilon_0}{2},
$$
by Corollary \ref{collapse}, $v_k'$ converges weakly in $W^{1,p}$
and $e^{v_k'}$ converges in $L^q$. Let $v'$ be the limit. We have
$$
-\Delta v'=e^{2v'},\s\int_{\R^2}e^{2v'}<+\infty.
$$  
By Theorem 1 in \cite{Chen-Li}, $v'=-\log(1+\frac{1}{4}|x-x_0|^2)$, and $\int_{\R^2}e^{2v'}=4\pi$, i.e. $(\C,e^{2v'}g_{euc})$ is a parametrization
of $(S^2\setminus \{y_0\},g_{S^2})$.

Put $\sigma_k=\pi^{-1}(x_k+r_k\pi(x))$, which defines a M$\ddot{o}$bius transformation of
$S^2$. If we set $g_k''=\sigma_k^*(g_k)=e^{2u_k''}g_{S^2}$, then 
$g_k''$ converges to $g_{S^2}$ weakly in $W^{1,p}_{loc}(S^2\setminus\{y_0\})$. 
Then we have
$$
\lim_{r\rightarrow 0}\lim_{k\rightarrow+\infty}\int_{B_r^{g_{S^2}}(y_0)}|K(g_k'')|d\mu_{g_k''}=0.
$$
Then $g_k''$ converges weakly in $W^{1,p}$.
\endproof

It is easy to deduce  Theorem \ref{main1} from Corollary \ref{main2} and the following lemma:

\begin{lem}
For any given 
$\Lambda$, if $\vol(S^2,g)\leq \Lambda$ , then for any $q>1$, there exit  $\tau>0$,  such that if 
$\|K(g)-1\|_{L^1(S^2,g)}<\tau$, then we can find a M$\ddot{o}$bius transformation $\sigma$, such that 
$\sigma^*(g)=e^{2u'}g_{S^2}$ with
$$
\|u'\|_{W^{1,p}}<C(p),\s \|e^{2u'}\|_{L^q}<C.
$$
\end{lem}

\subsection{A sequence with $K\geq 1$} Let $g_k$ be a metric sequence with $K(g_k)\geq 1$ and $\vol(g_k)\geq a>0$. We may
set $g_k=e^{2u_k}g_{S^2}$. We have
$$
diam(g_k)\leq\pi,\s \int_{S^2}|K(g_k)|d\mu_{g_k}=4\pi
$$ 
By the volume comparison theorem, we have $\vol(g_k)\leq 4\pi^2$.

Using the proof of Lemma \ref{sphere}, we may assume $u_k$ converges to $u$ weakly in $W^{1,p}
(S^2\setminus\mathcal{S})$, where $\mathcal{S}$ is a finite set.
Let $p\in\mathcal{S}$. If $\lim\limits_{r\rightarrow 0}\lim\limits_{k\rightarrow+\infty}\vol(B^{g_{S^2}}_r(p),g_k)>a$,
then we can find $x_k\rightarrow p$, $r_k\rightarrow 0$, and
$c_k$, such that
$u_k-c_k$ converges weakly in $W^{1,p}_{loc}(\C\setminus \mathcal{S'})$, where $\mathcal{S'}$ is a finite set. Then there exist a closed stable geodesic, which is 
impossible for $K\geq 1$.
Otherwise, we get $\lim\limits_{r\rightarrow 0}\lim\limits_{k\rightarrow+\infty}\vol(B^{g_{S^2}}_r(p),g_k)=0$.
By the volume comparison theorem, we have
$\lim\limits_{r\rightarrow 0}\lim\limits_{k\rightarrow+\infty} diam(B_r(p),g_k)\rightarrow 0$.
Then  we get  the following:
\begin{cor}
Let $g_k$ be a metric sequence with $K_{g_k}\geq 1$ and $\vol(g_k)\geq a>0$. The Gromov-Hausdorff limit of $(S^2,g_k)$ is a topological
sphere. Moreover, there exists a M$\ddot{o}$bius transformation 
$\sigma_k$, such that $\sigma_k^*(g_k)$ converges weakly
in $W^{1,p}_{loc}(S^2\setminus\mathcal{S},g_{S^2})$, where
$\mathcal{S}$ is a finite set.
\end{cor}

\end{document}